\documentclass[11pt]{article}
\usepackage{times,mathptm}
\usepackage{amsfonts,amscd,amssymb}
\newtheorem{lemma}{Lemma}[section]
\newtheorem{theorem}[lemma]{Theorem}
\newtheorem{proposition}[lemma]{Proposition}
\newtheorem{corollary}[lemma]{Corollary}

\newtheorem{remark}[lemma]{Remark}

\newtheorem{example}[lemma]{Example}

\newenvironment{proof}{{\it Proof.}}{\hfill $ \square $ \vskip 4mm}

\newcommand{\thlabel}[1]{\label{th:#1}}
\newcommand{\thref}[1]{Theorem~\ref{th:#1}}
\newcommand{\selabel}[1]{\label{se:#1}}
\newcommand{\seref}[1]{Section~\ref{se:#1}}
\newcommand{\lelabel}[1]{\label{le:#1}}
\newcommand{\leref}[1]{Lemma~\ref{le:#1}}
\newcommand{\prlabel}[1]{\label{pr:#1}}
\newcommand{\prref}[1]{Proposition~\ref{pr:#1}}

\newcommand{\relabel}[1]{\label{re:#1}}
\newcommand{\reref}[1]{Remark~\ref{re:#1}}
\newcommand{\exlabel}[1]{\label{ex:#1}}
\newcommand{\exref}[1]{Example~\ref{ex:#1}}

\newcommand{\eqlabel}[1]{\label{eq:#1}}
\newcommand{\eqref}[1]{(\ref{eq:#1})}

\def\bea{\begin{eqnarray*}}
\def\eea{\end{eqnarray*}}

\def\END{{\rm END}}
\def\INN{{\rm INN}}
\def\End{{\rm End}}
\def\Inn{{\rm Inn}}
\def\Gal{{\rm Gal}}
\def\Aut{{\rm Aut}}
\def\Hom{{\rm Hom}}
\def\Mod{{\rm Mod}}
\def\Pic{{\rm Pic}}
\def\PC{{\rm PC}}
\def\lan{\langle}
\def\ran{\rangle}

\def\units{{\mathbb G}_m}

\def\ol{\overline}

\def\ot{\otimes}

\def\phi{\varphi}
\def\ot{\otimes}

\input diagrams

\begin{document}
\title{On gradings of matrix algebras and descent theory}
\author {
S. Caenepeel\\ 
Faculty of Applied Sciences\\
Free University of Brussels, VUB\\ B-1050 Brussels, Belgium\and
S. D\u{a}sc\u{a}lescu and C. N\u{a}st\u{a}sescu\thanks{Research supported by the bilateral project
``Hopf Algebras in Algebra, Topology, Geometry and Physics" of the Flemish and
Romanian governments}\\
Faculty of Mathematics\\ University of Bucharest\\
RO-70109 Bucharest 1, Romania}

\date{}

\maketitle

\section{Introduction}
The aim of this paper is to investigate group gradings of matrix 
algebras. Let $k$ be a field and $G$ be a group. A $k$-algebra 
$A$ is called $G$-graded if $A=\oplus _{g\in G}A_g$, a direct sum of 
$k$-subspaces, such that $A_gA_h\subseteq A_{gh}$ for any $g,h\in G$. 
Several papers have been devoted to constructing and classifying 
group gradings on the matrix algebra ${\rm M}_m(k)$. One can view the 
algebra ${\rm M}_m(k)$ as a quotient of the path algebra of the quiver 
$\Gamma$, where $\Gamma$ is the complete graph on $m$ points. 
Certain gradings on ${\rm M}_m(k)$ were constructed 
in \cite{Green}, \cite{GreenMarcos} from weight functions on $\Gamma$. 
These gradings were called {\it good gradings} in \cite{dinr}, where they
were explained in an elementary way: they are 
the gradings for which all the matrix units $e_{ij}$ are homogeneous 
elements. The name good gradings indicates that it is easy to 
understand and construct them. As we will see in this paper, 
good gradings play 
a central role in the classification of all gradings on the matrix algebra. 
There is one more way to present good gradings. If $V$ is a 
$G$-graded vector space of dimension $m$, then there is a natural 
$G$-grading induced on $\End(V)$ (denoted by $\END(V)$), and this transfers 
to a grading of ${\rm M}_m(k)\cong \End(V)$. It turns out that this grading is 
a good grading, and moreover, any good grading can be obtained in this way.\\ 
In \seref{2} we prove that the gradings $\END(V)$ and $\END(W)$ are isomorphic 
if and only if the graded vector space $W$ is a suspension of $V$. This result 
can be used to classify all good $G$-gradings on ${\rm M}_m(k)$ by the orbits of 
the biaction of the symmetric group ${\cal S}_m$ and $G$ on the set $G^m$.  
We develop a combinatorial framework to use this for counting the good 
$G$-gradings in the case where $G=C_n$, the cyclic group of order $n$. \\
A fundamental result that we prove is that for an algebraically closed field 
$k$, any $C_n$-grading on ${\rm M}_m(k)$ is isomorphic to a good grading. We give two 
proofs for this. The first one relies on a result of \cite{dinr} concerned to graded 
Clifford theory, while the second one follows from an exact sequence of \cite{Beattie} 
concerned to graded Brauer groups. This fundamental result combined with  
the classification of the good $C_n$-gradings completes the classification of 
all $C_n$-gradings on ${\rm M}_m(k)$ for an algebraically closed field $k$. Our 
approach to the case where $k$ is arbitrary uses descent theory. 
Gradings of matrix algebras by the cyclic group $C_n$ have been investigated 
in \cite{bd} with a different approach under the assumption that the field 
$k$ contains a primitive $n$th root of unity. In this situation the 
group algebra $kC_n$ is a selfdual Hopf algebra, and then a $C_n$-grading 
is the same to an action of $C_n$. The method of \cite{bd} consists of 
using the Skolem-Noether theorem and the Jordan form for presenting  a 
$C_n$-grading. In this paper we do not need any assumption on the field. 
In \seref{1} we 
present some preliminary results about descent theory. We use them to prove 
the main result of \seref{3}, stating that for a field $k$ and a finite 
abelian group $G$, the forms of a certain (good) $G$-grading on 
${\rm M}_m(\overline{k})$ are in bijection to the isomorphism classes of Galois extensions 
of $k$ with a certain Galois group. We are able to give an explicit description of the 
grading corresponding to a Galois extension. We exemplify our theory by computing 
effectively the isomorphism types of $C_2$-gradings on ${\rm M}_2(k)$. This recovers and 
explains results of \cite{bd} in characteristic different from 2, and of 
\cite{Boboc} in characteristic 2. For proving the main result we need a graded version of
the  Rosenberg-Zelinsky exact sequence, and this has as a biproduct a graded version 
of the Skolem-Noether theorem.

\section{Descent theory}\selabel{1}
In \seref{3}, we will use classical descent theory, and we briefly
recall the basic results.
Our main references are \cite{KnusOjanguren} and \cite{Serre}.\\
Let $k$ be a commutative ring, and $l$ a commutative faithfully flat
$k$-algebra. Consider
an ``object" $N$ (a module, an algebra, a graded algebra,\ldots) defined
over $l$. Descent theory gives an answer to the following two problems:
\begin{enumerate}
\item{\textit{existence of forms}: does there exist an object $M$
defined over $k$ such that $M\ot_k l\cong N$?}
\item{\textit{classification}: if such a $k$-form of $N$ exists, classify all forms
up to isomorphism.}
\end{enumerate}
The answer to the first question is basically the following (cf. 
\cite[II.3.2]{KnusOjanguren}):
$M$ exists (as a $k$-module)
if and only if there exists a {\it descent datum} for $N$, this is
an $l\ot_k l$-isomorphism $u:\ l\ot_k N\to N\ot_k l$ such that $u_2=u_3\circ u_1$.
Here $u_1=I_l\ot u:\ l\ot_kl\ot_k N\to l\ot_kN\ot_k L$, and $u_2$ and $u_3$ are
defined in the same way.

If $N$ is an $l$-algebra and $u$ is an $l\ot_k l$-algebra map, 
then $M$ is a $k$ algebra. Similar properties hold for other structures like
$F$-graded algebras or $F$-module algebras ($F$ being a group).\\
In principle $M$ can be recovered from $N$, but this construction turns
out to be too complicated in practical computations. However, if $l$ is a
Galois extension of $k$ (in the sense of \cite{DI}), then we have a more
explicit construction. Write $G=\Gal(l/k)$, and take $\sigma\in G$.
A map $f\in \End_k(N)$ is called $\sigma$-semilinear if
$f(an)=\sigma(a)f(n)$ for all $a\in l$ and $n\in N$. With these notations,
we have the following result (cf. \cite[III.5.1]{KnusOjanguren}):

\begin{proposition}\prlabel{1.1} {\bf (Galois descent for modules)}
Let $l$ be a Galois extension of $k$ with group $G$. $N\in l$-mod descends
to $M\in k$-mod if and only if there exists a Galois descent datum for $N$,
that is a group homomorphism
$$\phi:\ G\to\Aut_k(N)$$
such that $\phi(\sigma)$ is $\sigma$-semilinear for all $\sigma\in G$. The
descended module $M$ is the fixed module $M=N^G$, with the $G$-action on $N$
induced by $\phi$.
\end{proposition}

Again, similar results hold for other structures. For example, if
$N=A$ is an $l$-algebra, then the descended module $M=B$ is a $k$-algebra
if all the $\phi(\sigma)$ are $k$-algebra maps.

Let us now look at the second problem. Assume that the $l$-module $N$ descends
to $M\in k$-mod. It is well-known (cf. \cite[II.8]{KnusOjanguren}) that we have a
bijection between the (pointed) set of isomorphism classes of $k$-modules
$M'$ that are isomorphic to $l\otimes_k M\cong N$ after base extension, and
the Amitsur cohomology set $H^1(l/k,\Aut_l(\bullet\ot_k M))$:
\begin{equation}\eqlabel{1.2.1}
{\rm Mod}(l/k,M)\cong H^1(l/k,\Aut_{\bullet}(\bullet\ot_k M))
\end{equation}
We refer to \cite{KnusOjanguren} for the definition of Amitsur cohomology. We obtain
similar statements about isomorphism classes of algebras, graded algebras,\ldots~
after we replace the automorphism group on the right hand side of
\eqref{1.2.1} by the appropriate automorphism group, for example, for
a $k$-algebra $A$,
the set of algebra forms of $l\otimes_k A$ becomes isomorphic to 
\begin{equation}\eqlabel{1.2.2}
{\rm Alg}(l/k,A)\cong H^1(l/k,\Aut_{\bullet}^{\rm Alg}(\bullet\ot_k A))
\end{equation}
Taking inductive limits over all faithfully flat extensions $l$ of $k$,
we obtain a formula for all twisted forms of $A$ as a $k$-algebra
\begin{equation}\eqlabel{1.2.3}
{\rm Alg}(k,A)\cong H^1(k_{\rm fl},\Aut_{\bullet}^{\rm Alg}(\bullet\ot_k A))
\end{equation}
Here are some well-known explicit variations of \eqref{1.2.3}, that will
be used in the sequel. Taking $M=k^n$ in \eqref{1.2.1}, we obtain
the (generalized) Hilbert Theorem 90, after
taking inductive limits:
\begin{equation}\eqlabel{1.2.4}
\Mod(l/k,k^n)\cong H^1(l/k,{\rm Gl}_n)
\end{equation}
If $k$ is a field, then both sides of \eqref{1.2.4} are trivial, since
there exists only one vector space of dimension $n$. For $n=1$, we find 
\begin{equation}\eqlabel{1.2.5}
\Pic(k)\cong H^1(l/k,\units)
\end{equation}
where, as usual, $\units(k)$ denotes the group of invertible elements of $k$.
Now let $F$ be a finite group. The dual $Fk=\oplus_{\gamma\in F} kv_{\gamma}$
of the group algebra $kF$ is an $F$-module algebra, the $F$-action is given
by
\begin{equation}\eqlabel{1.2.6}
\delta\cdot v_{\gamma}=v_{\gamma\delta^{-1}}
\end{equation}
for all $\gamma,\delta\in F$. $Fk$ is an $F$-Galois extension of $k$, and
an $F$-module algebra form of $Fk$ is again an $F$-Galois extension.
Furthermore
$$\Aut_F^{\rm alg}(Fk)=F(k)$$
the group of all continuous maps from ${\rm Spec}(k)$ to $F$, this is $F$
to the power the number of connected components of $k$. We find the
following formula for the set of commutative $F$-Galois extensions of $k$:
\begin{equation}\eqlabel{1.2.7}
\Gal(k,F)\cong H^1(k_{\rm fl},F)
\end{equation}

The correspondence \eqref{1.2.1} is known in principle, but the construction
is not very useful in practical situations. If $l/k$ is Galois, then we
have the following more explicit and practical description (see
\cite[Prop. 1.2]{CaenepeelDL} for the proof). 

\begin{proposition}\prlabel{1.3}
Let $l$ be a Galois extension of the commutative ring $k$, with Galois
group $G$, and consider a $k$-module $M$. The twisted form of $l\ot_k M$
corresponding to a cocycle $\psi\in Z^1(G,\Aut_l(l\ot_k M))$ under the
isomorphism
\begin{equation}\eqlabel{1.3.1}
\Mod(l/k,M)\cong H^1(G,\Aut_l(l\ot_k M))
\end{equation}
is equal to $(l\ot M)^G$, where the $G$-action on $l\ot_k M$ is induced
by the Galois descent datum $\phi:\ G\to \End_k(l\ot_k M)$, with 
\begin{equation}\eqlabel{1.3.2}
\phi(\sigma)(a\ot m)=\psi(\sigma)(\sigma(a)\ot m)=
\sigma(a)(\psi(\sigma)(1\ot m))
\end{equation}
for all $a\in l$ and $m\in M$. Similar results hold for the variations of
\eqref{1.3.1} (for algebras, coalgebras, Hopf algebras,\ldots).
\end{proposition}

\section{Good gradings}\selabel{2}
Let $k$ be a field, $m$ a positive integer and $G$ a group. We are 
interested to study $G$-gradings of the $k$-algebra 
of matrices ${\rm M}_m(k)$. A special class of gradings consists 
of the good gradings. A grading is called good if all the 
matrix units $e_{ij}$ are homogeneous elements. The class of good 
gradings is interesting since on one hand it is easy to  
construct them, and on the other hand any grading is isomorphic 
to a good one provided $k$, $G$ and $m$ satisfy certain conditions. 
This way to present good gradings from an interior point of view, by  
using certain elements of ${\rm M}_m(k)$, has the advantage that 
one can easily generate the good gradings by assigning arbitrary 
degrees to the elements $e_{12},e_{23},\ldots ,e_{m-1,m}$ 
(see \cite[Proposition 2.1]{dinr}). There is an alternative way to 
define good gradings, from an exterior point of view.\\
Let $V$ be a $G$-graded $k$-vector space, i.e. $V=\oplus _{g\in G}V_g$, 
a direct sum of vector subspaces of $V$. For any $\sigma \in G$ we 
define the left $\sigma$-suspension $V(\sigma )$ of $V$ as being the 
$G$-graded vector space which coincides with $V$ as a vector space, and 
having the grading given by $V(\sigma )_g=V_{g\sigma }$ for any 
$g\in G$. Similarly, the right $\sigma$-suspension $(\sigma )V$ of $V$ 
is just the space $V$ with the grading $(\sigma )V_g=V_{\sigma g}$ 
for any $g\in G$. Obviously $((\sigma )V)(\tau )=(\sigma )(V(\tau ))$ for 
any $\sigma ,\tau \in G$, and we simply denote this by $(\sigma )V(\tau )$.  
We consider the endomorphism algebra $\End(V)$ with map composition 
as multiplication. 
For any $\sigma \in G$ we denote by 
$$\END(V)_{\sigma }=\{ f\in \End(V)|f(V_g)\subseteq 
V_{\sigma g} \; {\rm for}\; {\rm any}\; g\in G\}$$ 
Then the sum $\sum _{\sigma \in G}\END(V)_{\sigma }$ is direct 
inside $\End(V)$, and we denote by 
$\END(V)=\oplus _{\sigma \in G}\END(V)_{\sigma}$, which is a $G$-graded 
algebra (see \cite{nv}). \\
If $V$ has dimension $m$, then 
$\END(V)=\End(V)\simeq {\rm M}_m(k)$, and this provides a $G$-grading 
on ${\rm M}_m(k)$ which turns out to be a good grading. In fact any 
good grading on ${\rm M}_m(k)$ is isomorphic to $\END(V)$ for some 
graded vector space $V$ of dimension $m$ (see \cite[Proposition 1.2]{dinr}), 
so from an exterior point of view we could define good gradings 
as gradings isomorphic to $\END(V)$ for some $V$. \\ 
Our first aim is to classify the good $G$-gradings.

\begin{theorem}  \label{clasend}
Let $V$ and $W$ be finite dimensional $G$-graded $k$-vector spaces. 
Then the graded algebras $END(V)$ and $END(W)$ are isomorphic 
if and only if there exists $\sigma \in G$ such that 
$W\simeq V(\sigma )$.
\end{theorem}

\begin{proof}
Let $v_1,\ldots ,v_m$ be a basis of $V$ consisting of 
homogeneous elements, say of degrees $g_1,\ldots ,g_m$. For any 
$1\leq i,j\leq m$ define $E_{ij}\in \END(V)$ by 
$E_{ij}(v_t)=\delta _{t,j}v_i$. Then $E_{ij}$ is a homogeneous 
element of degree $g_ig_j^{-1}$ of $\END(V)$,  
the set $(E_{ij})_{1\leq i,j\leq m}$ is a basis of $\END(V)$,  
and $E_{ij}E_{rs}=\delta_{j,r}E_{is}$ for any $i,j,r,s$. In particular 
$(E_{ii})_{1\leq i\leq m}$ is a complete system of orthogonal 
idempotents of $\END(V)$.\\
Let $u:\ \END(V)\rightarrow \END(W)$ be an isomorphism of $G$-graded algebras, 
and define $E'_{ij}=u(E_{ij})$ for any $1\leq i,j\leq m$. Denote 
$Q_i=Im(E'_{ii})$, which is a graded vector subspace of $W$. 
Since $(E'_{ii})_{1\leq i\leq m}$ is a complete system of orthogonal 
idempotents of $\END(W)$, we have that $W=\oplus _{1\leq i\leq m}Q_i$ and 
$E'_{ii}$ acts as identity on $Q_i$ for any $i$. Combined with the 
relation $E'_{ij}E'_{ji}=E'_{ii}$, this shows that 
$E'_{ij}$ induces an isomorphism of degree $g_ig_j^{-1}$ from 
$Q_i$ to $Q_j$. In particular $Q_j\simeq (g_1g_j^{-1})Q_1$ in  
$gr-k$ for any $j$. We obtain that $W\simeq \oplus _{1\leq j\leq m} 
(g_1g_j^{-1})Q_1$, showing that $Q_1$ has dimension 1. \\
Repeating this argument for the identity isomorphism from 
$\END(V)$ to $\END(V)$, and denoting by $R_1=Im(E_{11})$, we obtain that 
$V=\oplus _{1\leq j\leq m}(g_1g_j^{-1})R_1$. Since $Q_1$ and $R_1$ are 
graded vector spaces of dimension 1, we have that $Q_1\simeq R_1(\sigma )$ 
for some $\sigma \in G$. Hence 
$W\simeq \oplus _{1\leq j\leq m}(g_1g_j^{-1})Q_1
\simeq \oplus _{1\leq j\leq m}(g_1g_j^{-1})(R_1(\sigma )) 
\simeq V(\sigma )$.\\
For the other way around, it is easy to see that 
$\END(V)\simeq \END(V(\sigma ))$ as $G$-graded algebras.
\end{proof} 

The previous theorem can be used to classify all good $G$-gradings 
on the matrix algebra ${\rm M}_m(k)$. Indeed, any such grading is of the 
form $\END(V)$ for some $G$-graded vector space $V$ of dimension $m$. 
To such a $V$ we can associate an $m$-tuple 
$(g_1,\ldots ,g_m)\in G^m$ consisting of the degrees of the 
elements in a homogeneous basis of $V$. Conversely, to any such 
a $m$-tuple, one can associate a $G$-graded vector space 
of dimension $m$. Obviously, the $G$-graded vector space 
associated to a $m$-tuple coincides to the one associated 
to a permutation of the $m$-tuple. In fact Theorem \ref{clasend} 
can be reformulated in terms of $m$-tuples as follows.  
If $V$ and $W$ are $G$-graded vector spaces of dimension $m$ associated 
to the $m$-tuples $(g_1,\ldots ,g_m)$ and $(h_1,\ldots ,h_m)$, 
then $\END(V)\simeq \END(W)$ as $G$-graded algebras if and only if 
there exist $\sigma \in G$ and $\pi $ a permutation of $\{ 1,\ldots ,m\}$ 
such that $h_i=g_{\pi (i)}\sigma$ for any $1\leq i\leq m$. 
We have showed the following.

\begin{corollary}  \label{clasbiaction}
The good $G$-gradings of ${\rm M}_m(k)$ are classified by the orbits of 
the biaction of the symmetric group ${\cal S}_m$ (from the left) and 
$G$ (by translation from the right) on the set $G^m$.
\end{corollary}

\subsubsection*{Classification of good gradings over a cyclic group}
In the particular case where $G=C_n=<c>$, a cyclic group with $n$ elements, 
we are able to count the orbits of this biaction. Let 
$F_{n,m}$ be the set of all $n$-tuples $(k_0,k_1,\ldots ,k_{n-1})$ of 
non-negative integers with the property that $k_0+k_1+\ldots +k_{n-1}=m$.
It is a wellknown combinatorial fact that $F_{n,m}$ has 
${m+n-1}\choose{n-1}$ elements. To any element $(g_1,\ldots ,g_m)\in G^m$  
we can associate an element $(k_0,k_1,\ldots ,k_{n-1})\in F_{n,m}$ such 
that $k_i$ is the number of appearances of the element $c^i$ in 
the $m$-tuple $(g_1,\ldots ,g_m)$ for any $1\leq i\leq n$. By  
Corollary \ref{clasbiaction} we see that the number of orbits of the 
$({\cal S}_m,G)$-biaction on $G^m$ is exactly the number of orbits 
of the left action by permutations of the subgroup $H=<\tau >$ of 
${\cal S}_n$ on the set $F_{n,m}$, where $\tau$ is the cyclic permutation  
$(1\; 2\; \ldots n)$. This number of orbits is the one we will 
effectively compute. \\
If $\alpha =(k_0,k_1,\ldots ,k_{n-1})\in F_{n,m}$ 
and $d$ is a positive divisor of $n$, then $\tau ^d\alpha =\alpha$  
if and only if $k_{i+d}=k_i$ for any $0\leq i\leq n-1-d$, i.e. the 
first $d$ positions of $\alpha$ repeat $\frac{n}{d}$ times. 
In particular we must 
have that $\frac{n}{d}$ divides $m$. For any such $d$, let us denote 
by $A_d$ the set of all $\alpha \in F_{n,m}$ stabilized by $\tau ^d$. 
Since for $\alpha \in A_d$ we have that $k_0+k_1+\ldots +k_{d-1}= 
\frac{md}{n}$, we see that $A_d$ has 
${\frac{md}{n}+d-1}\choose{d-1}$ elements.\\
Let ${\cal D}(n,m)$ be the set of all positive divisors $d$ of $n$ 
with the property that $n/d$ divides $m$. If $d_1,d_2\in {\cal D}(n,m)$, 
then $(d_1,d_2)\in {\cal D}(n,m)$. Indeed, since $n$ divides $d_1m$ and 
$d_2m$, then $n$ also divides $(d_1m,d_2m)=(d_1,d_2)m$, therefore 
$(d_1,d_2)\in {\cal D}(n,m)$. It follows that ${\cal D}(n,m)$ is a lattice 
with the order given by divisibility. 
For any $d\in {\cal D}(n,m)$ we denote by ${\cal D}(n,m,d)$ the set 
of all elements $d'$ of ${\cal D}(n,m)$ which divide $d$, and by  
${\cal D}_0(n,m,d)$ the set of all maximal elements of ${\cal D}(n,m,d)$. 
The following is immediate. 

\begin{lemma} \label{lemaint}
For any $d_1,d_2\in {\cal D}(n,m)$ we have that $A_{d_1}\cap A_{d_2}=
A_{(d_1,d_2)}$.
\end{lemma}

The following describes the elements with the orbit of length $d$. 

\begin{lemma} \label{lengthd}
For any $d\in {\cal D}(n,m)$ denote by $B_d$ the set of all the 
elements of $F_{n,m}$ having the orbit of length $d$. Then 
$$B_d=A_d-\bigcup _{d'\in {\cal D}(n,m,d)}A_{d'}=
A_d-\bigcup _{d\in {\cal D}_0(n,m,d)}A_{d'}$$
\end{lemma}

\begin{proof}
The orbit of an element $\alpha$ has length $d$ 
if and only if the stabilizer of $\alpha$ is a subgroup 
with $\frac{n}{d}$ elements of $H$, thus equal to 
$<\tau ^d>$. The result follows now from the definition of $A_d$.
\end{proof} 

\begin{corollary}
Let $d\in {\cal D}(n,m)$ and $p_1,\ldots ,p_s$ all the distinct 
prime divisors of $d$ such that $\frac{d}{p_1},\ldots ,\frac{d}{p_s}
\in {\cal D}(n,m)$. Then 
$$|B_d|=|A_d|+\sum _{1\leq t\leq s}\sum_{1\leq i_1<\ldots <i_t\leq s}
(-1)^t|A_{\frac{d}{p_{i_1}\ldots p_{i_t}}}|$$
and this is known taking into account the fact that for any $d'$ we have 
that 
$|A_{d'}|={{\frac{md'}{n}+d'-1}\choose{d'-1}}$.
\end{corollary}

\begin{proof}
We have that ${\cal D}_0(n,m,d)=\{ \frac{d}{p_1},\ldots ,
\frac{d}{p_s}\}$. The result follows now from Lemma \ref{lengthd} and 
by applying the principle of inclusion and exclusion. 
\end{proof}

\begin{theorem}
The number of isomorphism types of good $C_n$-gradings of the 
algebra ${\rm M}_m(k)$ is 
$$\sum _{d\in {\cal D}(n,m)}\frac{1}{d}|B_d|
=\sum _{d\in {\cal D}(n,m)}\frac{1}{d}(|A_d|-
|\bigcup _{d\in {\cal D}_0(n,m,d)}A_{d'}|)$$
\end{theorem}

\begin{proof} The number of isomorphism types of the good gradings 
is the number of orbits of the action of $H$ on $F_{n,m}$. 
We have seen that if the orbit of an element $\alpha \in 
F_{n,m}$ has length $d$, then necessarily $d\in {\cal D}(n,m)$. 
The result follows since the number of orbits of length $d$ is 
$\frac{1}{d}|B_d|$.
\end{proof}

\begin{example}\rm
Let $n=p^r$ with prime $p$ and let $m$ be a positive integer. 
We define $q$ by $q=r$ in the case where $n$ divides $m$,  
or  $q$ is the exponent of $p$ in $m$ in the case where $n$  
does not divide $m$. Then 
${\cal D}(n,m)=\{ p^i|r-q\leq i\leq r\}$ and 
the number of isomorphism types 
of good $C_n$-gradings on ${\rm M}_m(k)$ is 
$$\frac{1}{p^{r-q}}{{\frac{m}{p^{q}}+p^{r-q}-1}\choose{p^{r-q}-1}}+
\sum _{r-q< i\leq r}\frac{1}{p^i}
({{\frac{m}{p^{r-i}}+p^i-1}\choose{p^i-1}}-
{{\frac{m}{p^{r-i-1}}+p^{i-1}-1}\choose{p^{i-1}-1}})
$$
In particular, if $n=p$, then the number of isomorphism types 
of good $C_p$-gradings on ${\rm M}_m(k)$ is 
$1+\frac{1}{p}({{m+p-1}\choose{p-1}}-1)$ if $p$ divides $m$, 
and $\frac{1}{p}{{m+p-1}\choose{p-1}}$ if $p$ does not divide 
$m$. This fact was proved in \cite[Proposition 3.3]{bd} in the case where 
$k$ contains a primitive $p$-th root of $1$.
\end{example}

\subsubsection*{The algebraically closed field case}
We show that in certain cases all gradings of the matrix algebras 
are isomorphic to good ones.

\begin{theorem}\thlabel{2.8}
Let $k$ be an algebraically closed field and $m,n$ positive integers. 
Then any $C_n$-grading of the matrix algebra ${\rm M}_m(k)$ is isomorphic to a 
good grading.
\end{theorem}

\begin{proof} Let $R={\rm M}_m(k)$ be the matrix algebra endowed with a 
certain $C_n$-grading, and pick $\Sigma \in R-gr$ a graded simple 
module. Then $\Delta =End_R(\Sigma )=END_R(\Sigma )$ is a $C_n$-graded 
algebra. Moreover, $\Delta$ is a crossed product when regarded as 
an algebra graded by the support of the $C_n$-grading of $\Delta$. 
Thus $\Delta \simeq \Delta _1\sharp _{\sigma}H$ for a cyclic group 
$H$ and a cocycle $\sigma$. On the other hand $\Delta _1=End_{R-gr}
(\Sigma )$ is a finite field extension of $k$, so $\Delta _1=k$. 
Hence $\Delta$ is a crossed product of $\Delta _1$, which is central 
in $\Delta$, and the cyclic group $H$, so $\Delta$ is commutative.\\
Since $R$ is a semisimple algebra with precisely one isomorphism type of 
simple module, say $S$, we have that $\Sigma \simeq S^p$ as 
$R$-modules for some positive integer $p$. But then 
$\Delta \simeq End_R(S^p)\simeq {\rm M}_p(k)$, and the commutativity of 
$\Delta$ shows that $p$ must be $1$. Then $\Sigma \simeq S$, so 
there exists a graded $R$-module which is simple as an $R$-module.  
By \cite[Theorem 1.4]{dinr}, the grading is isomorphic to a good  
one.
\end{proof}

\section{Non-good gradings and descent theory}\selabel{3}
We have seen that all $C_n$-gradings on matrix algebras over an
algebraically closed field are good.
In this Section, we will classify $C_n$-gradings on matrix algebras
over an arbitrary field, using descent. More precisely, we will
discuss the following problem: given a finite abelian group $G$, 
and a $G$-graded $k$-vector space $V$, describe all gradings on 
$A=\End_k(V)$ that become isomorphic to $\END_{\ol{k}}(\ol{V})$ after a base extension,
i.e. all gradings that are forms of the good grading 
$\END_{\ol{k}}(\ol{V})$, where $\overline{k}$ stands for the algebraic closure of
$k$ and $\overline{V}=\overline{k}\otimes _kV$.

First, we introduce the following notation for a $G$-graded
vector space $V$:
\begin{equation}\eqlabel{3.1.1}
{\cal I}(V) =\{\sigma\in G~|~V(\sigma)\cong V~{\rm as~graded~vector~spaces}\}
\end{equation}
It is obvious that ${\cal I}(V)$ is a subgroup of $G$. Also
remark that $\sigma\in {\cal I}(V)$ if and only if there exists
an automorphism $f$ of $V$ that is homogenous of degree $\sigma$.

\begin{lemma}\lelabel{3.1}
Let $k$ be a field, $G$ a finite abelian group, and $V$ a $G$-graded
vector space with ${\cal I}(V)=H$. Then there exists a $G$-graded
vector space $W$ such that
$$V\cong kH\ot W~~{\rm and}~~{\cal I}(W)=\{e\}$$
\end{lemma}

\begin{proof}
We choose coset representatives $e=\sigma_1,\sigma_2,\cdots, \sigma_k$
of $H$ in $G$, and let
$$W=\oplus_{i=1}^k V_{\sigma_i}$$
For all $\sigma\in H$, we can find an automorphism $f_{\sigma}$ of $V$
of degree $\sigma$, and this implies that 
$V_{\sigma_i}$ and $V_{\sigma\sigma_i}$ are isomorphic, and have the
same dimension. This can also be stated in the following way:
$$\oplus_{\sigma\in H\sigma_i}V_{\sigma}\cong kH\ot V_{\sigma_i}$$
and it follows easily that $V\cong kH\ot W$.\\
Now take $\sigma\in {\cal I}(W)$. There exists an isomorphism $f:\ W\to 
W$ of degree $\sigma$, and $I_{kH}\ot f:\ V\cong kH\ot W\to V\cong kH\ot W$
is also an isomorphism of degree $\sigma$, and this implies that
$\sigma\in H$. If $\sigma\neq e$, then $f$ sends $V_{\sigma_i}$ to some
$V_{\sigma_j}$, with $i\neq j$, implying that $\sigma_j=\sigma_i\sigma$,
and contradicting $\sigma\neq e$. As a consequence, $\sigma=e$, and
${\cal I}(W)=\{e\}$.
\end{proof}

Now we fix a $G$-graded vector space $V$, and write $A=\END_k(V)$.
$\Aut_{\rm gr}(A)$ will denote the group of algebra
automorphisms of degree $e$ of $A$, and $\INN(A)$ the subgroup consisting of
inner algebra automorphisms induced by an isomorphism
$f$ of degree $e$ of $V$. $\PC(k,G)$ will be the group of grade preserving isomorphism
classes of dimension one graded vector spaces. It is well-known and easy to
show that $\PC(k,G)\cong G$: the graded vector space $k(\sigma)$ represents the
isomorphism class in $\PC(k,G)$ corresponding to $\sigma\in G$.

The {\sl Skolem-Noether Theorem} tells us that every automorphism of a matrix
ring is inner. If we work over a commutative ring instead of a field,
then this is no longer true, and the failure of the Skolem-Noether Theorem
is measured by the {\sl Picard group} of the ring, which is the set of
isomorphism classes of rank one projective modules. In fact we have
an exact sequence, known as the {\sl Rosenberg-Zelinsky exact sequence}
\begin{equation}\eqlabel{3.2.1}
0\rTo^{}\Inn(A)\rTo^{}\Aut(A)\rTo^{\phi}\Pic(k)
\end{equation}
with $\phi(\gamma)$ represented by $I=\{x\in A~|~x\gamma(a)=ax~
{\rm for~all~}a\in A\}$
and ${\rm Im}(\phi)=\{[I]\in\Pic(k)~|~I\ot V\cong V\}$ (see
\cite[IV.1.3]{KnusOjanguren} for the proof). Allthough we are interested in
algebras over a field, we mention the commutative ring case because the
graded version of the Picard group is not trivial, even if we work over
a field: it is exactly the $\PC(k,G)$ introduced above. We will need the
following graded version of the Rosenberg-Zelinsky exact sequence
(cf. \cite[13.6.1]{Caenepeel}):

\begin{proposition}\prlabel{3.2}
Let $G$ be a finite abelian group, $V$ a $G$-graded vector space over
a field $k$, and $A=\END_k(V)$. Then we have an exact sequence
$$0\rTo^{}\INN(A)\rTo^{}\Aut_{\rm gr}(A)\rTo^{\phi}_{}{\cal I}(V)\rTo^{}0$$
\end{proposition}

\begin{proof}
A straightforward adaption of the proof exhibited in \cite[IV.1.3]{KnusOjanguren}
gives an exact sequence
$$0\rTo^{}\INN(A)\rTo^{}\Aut_{\rm gr}(A)\rTo^{\phi}\PC(k,G)$$
with
$${\rm Im}(\phi)=\{[I]\in\PC(k,G)~|~I\ot V\cong V\}$$
Now we have seen that $\PC(k,G)\cong G$, and it is clear that
$k(\sigma)\ot V\cong V$ if and only if $\sigma\in {\cal I}(V)$,
and it follows that ${\rm Im}(\phi)\cong {\cal I}(V)$.
\end{proof}

\begin{remark}\rm\relabel{3.3}
The map $\phi:\ \Aut_{\rm gr}(A)\to \PC(k,G)$ can be described
as follows: $\phi(\gamma)=[I]$, with
$$I=\{x\in A~|~x\gamma(a)=ax~
{\rm for~all~}a\in A\}$$
Identifying $\PC(k,G)$ and $G$, this means that $\phi(\gamma)=\sigma$
if and only if $\gamma$ is induced by an isomorphism of degree $\sigma$.
This leads to the graded version of the Skolem-Noether Theorem: 
every isomorphism of a good graded matrix ring is induced by a homogeneous
matrix.
\end{remark}

Our next step is to introduce the category ${\cal G}_k$
of $G$-gradings on matrix
rings over a field $k$. The objects are pairs
$$(V,\{A_{\sigma}~|~\sigma\in G\})$$
where $V$ is a finite dimensional $k$-vector space, and
$$A=\End_k(V)=\oplus_{\sigma\in G}A_{\sigma}$$
is a $G$-graded $k$-algebra. An (iso)morphism 
$$(V,\{A_{\sigma}~|~\sigma\in G\})\to (V',\{A'_{\sigma}~|~\sigma\in G\})$$
consists of an isomorphism of vector spaces $f:\ V\to V'$
such that the induced map $\End(f):\ A\to A'$ is a
$k$-algebra isomorphism of degree $e$. The good gradings
form a full subcategory ${\cal GG}_k$ of ${\cal G}_k$.
In the situation where $G=C_n$ and $k$ is algebraically closed,
all gradings are good, and the two categories coincide.
We now want to apply \eqref{1.2.3}, but with the category of
$k$-algebras replaced by the category ${\cal G}_k$.
First we have to compute the automorphisms of a good grading
in the category ${\cal G}_k$.\\
Let $V$ be a graded vector space, and $f:\ V\to V$
an isomorphism of the good grading $(V,\END_k(V))\in {\cal GG}_k$.
In \reref{3.3}, we have seen that $f$ is homogeneous of degree
$\sigma$, for some $\sigma\in {\cal I}(V)\subset G$.
More precisely, we
have the following commutative diagram:
\begin{equation}\eqlabel{3.2.2}
\begin{diagram}
&&0&&&&&&\\
&&\dTo^{}&&&&&&\\
&&\units(k)&&0&&&&\\
&&\dTo^{}&&\dTo^{}&&&&\\
0&\rTo^{}&\Aut_{\rm gr}(V)&\rTo^{}&\Aut_{{\cal G}_k}(V,\END_k(V))&\rTo^{}&
{\cal I}(V)&\rTo^{}&0\\
&&\dTo^{}&&\dTo^{}&&\dTo^{=}&&\\
0&\rTo^{}&\INN(\END_k(V))&\rTo^{}&\Aut_{\rm gr}(\END_k(V)&\rTo^{\phi}&
{\cal I}(V)&\rTo^{}&0\\
&&\dTo^{}&&&&&&\\
&&{0}&&&&&&
\end{diagram}
\end{equation}
It follows from \prref{3.2} and some diagram chasing that the
rows and columns in the diagram are exact. Now we know from
\eqref{1.2.3} (applied to the category of gradings), that the gradings
on $\End_k(V)$ that become isomorphic to the good grading
$\END_k(V)$ after a base extension are described by the first
cohomology group $H^1(\ol{k}/k,\Aut_{{\cal G}_{\ol k}}(\ol{V},\END_{\ol{k}}(\ol{V}))$,
with $\ol{V}=\ol{k}\ot V$. From the exactness of the second
row in \eqref{3.2.2}, it follows that we have an exact sequence of
cohomology
$$H^1(\ol{k}/k,\Aut_{\rm gr}(\ol{V}))\to
H^1(\ol{k}/k,\Aut_{{\cal G}_{\ol k}}(\ol{V},\END_{\ol{k}}(\ol{V})))
\rTo^{\alpha}H^1(\ol{k}/k,{\cal I}(\ol{V}))$$
Now we can write
$$\ol{V}=\oplus_{\sigma\in G}(k^{n(\sigma)})(\sigma)$$
where $n(\sigma)$ is the dimension of the part of degree $\sigma$ of $V$;
this implies that
$$\Aut_{\rm gr}(V)=\oplus_{\sigma\in G}{\rm Gl}_{n(\sigma)}(\ol{k})$$
and $H^1(\ol{k}/k,\Aut_{\rm gr}(\ol{V}))$ becomes trivial, by the
generalized Hilbert 90, cf. \eqref{1.2.4}. This implies that
the map $\alpha$ is injective. In fact we have the following:

\begin{theorem}\thlabel{3.4}
Let $k$ be a field, $G$ a finite abelian group, and $V$ a $G$-graded
$k$-vector space. Then the map $\alpha$ described above is an isomorphism.
Consequently, we have a bijective correspondence between
\begin{itemize}
\item isomorphism classes of gradings on $\End_k(V)$ that become
isomorphic to the grading $\END_k(V)$ after base extension;
\item isomorphism classes of
Galois extensions $l$ of $k$ with Galois group $H={\cal I}(V)$.
\end{itemize}
\end{theorem}

\begin{proof}
We have already seen that $\alpha$ is injective, so it suffices to show
that $\alpha$ is surjective. The second statement then follows from
the fact that the isomorphism classes of gradings are in bijective
correspondence with  
$H^1(\ol{k}/k,\Aut_{{\cal G}_{\ol k}}(\ol{V},\END_{\ol{k}}(\ol{V}))$,
and isomorphism classes of $H$-Galois extensions are in bijective correspondence
with $H^1(\ol{k}/k,H)$, cf. \eqref{1.2.7}.\\
Consider a cocycle in $Z^1(\ol{k}/k,H)$. This cocycle corresponds to
an $H$-Galois extension $l$ of $k$, and $l$ splits this cocycle,
i.e. it is represented by some cocycle in
$$Z^1(H=\Gal(l/k),H)=\Hom(H,H)$$
and this cocycle is exactly the identity map $I_H:\ H\to H$. To prove
that $\alpha$ is surjective, it suffices to find an inverse image
in $H^1(H,\Aut_{{\cal G}_l}(l\ot V,\END_{l}(l\ot V)))$. Equivalently,
we may look for the corresponding descent datum (see \prref{1.3})
$$\phi:\ H\to \Aut_{{\cal G}_l}(l\ot V,\END_{l}(l\ot V))$$
$H$ is a finite abelian group, so we can write
$$H=\lan \sigma_1\ran \times\cdots\times \lan \sigma_n\ran$$
with $\sigma_i^{m_i}=e$. We know $\phi$ once we know all the
$\phi(\sigma_i)$. Recall from \leref{3.1} that
$$V\cong kH\ot W$$
so
$$l\ot V\cong l\ot k\lan\sigma_1\ran\ot\cdots\ot k\lan\sigma_n\ran\ot W$$
We define $\phi(\sigma_i):\ l\ot V\to l\ot V$ as follows:
$$\phi(\sigma_i)=\sigma_i\ot I_{k\lan \sigma_1\ran}\ot\cdots\ot \varphi_i
\ot\cdot\ot I_{k\lan \sigma_n\ran}\ot I_W$$
with $\varphi_i:\ k\lan \sigma_i\ran\to k\lan \sigma_i\ran$
given by
$$\varphi_i(\sigma_i^l)=\sigma_i^{l+1}$$
i.e. the matrix of $\varphi_i$ with respect to the basis
$\{e,\sigma_i,\sigma_i^2,\cdots,\sigma_i^{m_i-1}\}$ is
$$\pmatrix{
0&1&0&\cdots&0&0\cr
0&0&1&\cdots&0&0\cr
\vdots&\vdots&\vdots&&\vdots&\vdots\cr
0&0&0&\cdots&0&1\cr
1&0&0&\cdots&0&0\cr}$$
Clearly $\phi(\sigma_i)$ is $\sigma_i$-semilinear, and $\phi(\sigma_i)^{m_i}=1$,
which means that $\phi$ is a descent datum. For all $\sigma\in H$,
$\phi(\sigma)$ is homogeneous of degree $\sigma$, and the same thing holds
for the corresponding cocycle $\psi$. This means exactly that $\alpha$ maps $\psi$
to the identity map $I_H$, proving that $\alpha$ is surjective.
\end{proof}

\begin{remark}\rm\relabel{3.5}
Let $l$ be an $H$-Galois extension of $k$. Let us give an explicit
description of the grading corresponding to $l$ on $\End_k(kH)$.
To this end, we first recall the definition of a Galois extension
of a commutative ring (see \cite[III.1.2]{DI}). Let $l/k$ be an extension
of commutative rings, and assume that $H$ acts on $l$ as a group
of $k$-automorphisms. Now we consider the graded $k$-algebra
$$\Delta(l,H)=\oplus_{\sigma\in H} lu_{\sigma}$$
with multiplication
$$(au_{\sigma})(bu_{\tau})=a\sigma(b)u_{\sigma\tau}$$
and $H$-grading
$${\rm deg}(u_{\sigma})=\sigma$$
We have a $k$-algebra map
$$j:\ \Delta(l,H)\to\End_k(l)~~;~~j(au_{\sigma})(b)=a\sigma(b)$$
and $l$ is a Galois extension of $k$ with group $H$ if and only if
$l^H=k$ and $j$ is an isomorphism. If $j$ is an isomorphism, then the
grading on $H$ induces a grading on $\End_k(l)$, and this is the
grading we are looking for.\\
Let us show that this grading is a form of the good grading on
$\End_k(kH)$. Since $l$ is a Galois extension of $k$ with group
$H$, we have an isomorphism
$$m:\ l\ot l\to Hl~~;~~m(a\ot b)=\sum_{\sigma\in H}a\sigma(b)v_{\sigma}$$
(see \cite[III.1.2]{DI}). This means in fact that $l$ is a form of the
dual of the groupalgebra $Hk$, and it suffices to show that we have
a graded isomorphism
$$f:\ \Delta(Hk,H)\to \End_k(kH)$$
Let $\{e_{\sigma,\tau}~|~\sigma,\tau\in H\}$ be the canonical basis
of $\End_k(kH)$ consisting of elementary matrices with respect to the
basis $H$ of $kH$, i.e.
$$e_{\sigma,\tau}e_{\rho,\nu}=\delta_{\tau,\rho}e_{\sigma,\nu}~~{\rm and}~~
{\rm deg}(e_{\sigma,\tau})=\tau\sigma^{-1}$$
$\Delta(Hk,H)$ has $k$-basis $\{v_{\sigma}u_{\tau}~|~\sigma,\tau\in H\}$,
with
$$(v_{\sigma}u_{\tau})(v_{\rho}u_{\nu})=\delta_{\sigma,\rho\tau^{-1}}
v_{\sigma}u_{\tau\nu}~~{\rm and}~~{\rm deg}(v_{\sigma}u_{\tau})=\tau$$
Now define
$$f(e_{\sigma,\tau})=v_{\sigma}u_{\tau\sigma^{-1}}$$
It is straightforward to show that $f$ is a grade preserving isomorphism
of $k$-algebras.
\end{remark}

\begin{remark}\rm\relabel{3.6}
Our results are closely related to an exact sequence of Beattie
\cite{Beattie}, which we will recall here in the particular situation
that is of interest to us. Let $k$ be a field, $G$ a finite abelian
group, and look at $G$-graded central simple algebras (in particular
graded matrix algebras). Two such graded algebras $A$ and $B$ are called graded
equivalent if there exist matrix rings ${\rm M}_n(k)$ and ${\rm M}_m(k)$
equipped with a good grading such that $A\ot {\rm M}_n(k)$ and
$B\ot {\rm M}_m(k)$ are isomorphic as graded algebras. Equivalence classes
of graded central simple algebras form a group under the operation
induced by the tensor product, which we call the $G$-graded Brauer
group of $k$, and we denote it by ${\rm BC}(k,G)$. Beattie \cite{Beattie}
proved that we have a split exact sequence
\begin{equation}\eqlabel{3.6.1}
1~\to~{\rm Br}(k)~\to~{\rm BC}(k,G)~\to~{\rm Gal}'(k,G)~\to~1
\end{equation}
Here ${\rm Gal}'(k,G)$ consists of all isomorphism of (possibly
noncommutative) Galois extensions of $k$ with group $G$. The
previously introduced ${\rm Gal}(k,G)$ contains only commutative
Galois extensions. If $G$ is cyclic, then all Galois extensions
are commutative. If $k=\ol{k}$ is algebraically closed, then
${\rm Gal}(k,G)$ is trivial. In particular, it follows
that ${\rm Br}(\ol{k})={\rm BC}(\ol{k},C_n)$, implying that any
$C_n$-grading on a matrix algebra over an algebraically closed
field is good, providing an alternative proof of \thref{2.8}.
The inverse image of a (commutative) Galois extension $l$ of
$k$ in \eqref{3.6.1} is represented by the graded algebra constructed
in \reref{3.5}.
\end{remark}

\begin{example}\exlabel{3.7}\rm
In some cases, we can do the descent of the grading explicitely.
Let $k$ be a field of characteristic different from $2$, and take
$G=C_2=\{e,\sigma\}$. We will look at forms of the nontrivial good grading on
${\rm M}_2(k)$, given by
$${\rm deg}(e_{11})={\rm deg}(e_{22})=e~~;~~
{\rm deg}(e_{21})={\rm deg}(e_{12})=\sigma$$
which is also the good grading $\END_k(V)$, with
$V=ke\oplus k\sigma$ with the obvious grading. 
From \thref{3.4}, we
know that isomorphism classes of these gradings are given by the
quadratic extensions
${\rm Gal}(k,C_2)\cong k^*/(k^*)^2$. Recall that the quadratic extension
corresponding to $\alpha\in k^*$ is
$$l=k[X]/(X^2-\alpha)$$
Let $\sqrt{\alpha}$ be the element in $l$ represented by $X$. The descent
datum corresponding to $l$ is
$$\sigma\ot\varphi:\ l\ot V\to l\ot V$$
where the matrix of $\varphi$ is
$$\pmatrix{0&1\cr 1&0\cr}$$
The induced map $\phi=\End(\varphi):\ {\rm M}_2(k)\to {\rm M}_2(k)$ is given by
$$\phi\pmatrix{a&b\cr c&d\cr}=\pmatrix{d&c\cr b&a\cr}$$
The descended graded matrix ring consists of the matrices in
${\rm M}_2(l)$ that are fixed by $\sigma\ot\phi$, it is
$$A=\{\pmatrix{a&b\cr \sigma(a)&\sigma(b)\cr}~|~a,b\in l\}$$
or
$$A=\{\pmatrix{a+b\sqrt{\alpha}&c+d\sqrt{\alpha}\cr
c-d\sqrt{\alpha}& a-b\sqrt{\alpha}\cr}~|~a,b,c,d\in l\}$$
$A$ is the $k$-endomorphism ring of the descended module $\tilde{V}$,
which consists of the elements in $l\ot V$ fixed by $\sigma\ot\varphi$:
$$\tilde{V}=\{(a+b\sqrt{\alpha})e+(a-b\sqrt{\alpha})\sigma~|~a,b\in k\}$$
A $k$-basis for $\tilde{V}$ is $\{u=e+\sigma,v=e-\sigma\}$.\\
Now we observe that
$$A_e=k\pmatrix{1&0\cr 0&1\cr}\oplus k\sqrt{\alpha}\pmatrix{1&0\cr 0&-1\cr}~~;~~
A_{\sigma}=k\pmatrix{0&1\cr 1&0\cr}\oplus k\sqrt{\alpha}\pmatrix{0&1\cr 1&0\cr}$$
Switching to the new basis $\{u,v\}$, we find the following grading
on $B={\rm M}_2(k)$:
$$B_e=k\pmatrix{1&0\cr 0&1\cr}\oplus k\pmatrix{0&1\cr \alpha &0\cr}~~;~~
B_\sigma=k\pmatrix{1&0\cr 0&-1\cr}\oplus k\pmatrix{0&1\cr -\alpha &0\cr}$$
and we have recovered \cite[Corollary 3.6]{bd}
\end{example}

\begin{example}\exlabel{3.8}\rm
Now let $k$ be a field of characteristic $2$. Again we describe the grading
forms on ${\rm M}_2(k)$ with the same good grading as in \exref{3.7}. The difference
with \exref{3.7} is that the description of the isomorphism classes of
quadratic extension is different. We now have
$${\rm Gal}(k,C_2)\cong k/\{x^2-x~|~x\in k\}$$
The quadratic extension corresponding to $\alpha\in k$ is now
$$l=k[X]/(X^2-X-\alpha)$$
We write $x$ for the element in $l$ represented by $X$. The action of
$C_2$ on $l$ is given
$$\sigma(x)=x+1$$
Proceeding as
in \exref{3.7}, we now find for the descended graded matrix ring
$$A=\{\pmatrix{a+bx& c+dx\cr c+d(x+1)& a+b(x+1)\cr}~|~a,b,c,d\in k\}$$
The grading is now described by the formulas
$$A_e=k\pmatrix{1&0\cr 0&1\cr}\oplus k\pmatrix{x&0\cr 0&x+1\cr}~~;~~
A_{\sigma}=k\pmatrix{0&1\cr 1&0\cr}\oplus k\pmatrix{0&x\cr x+1&0\cr}$$
and $A$ is the $k$-endomorphism ring of the $k$-module
$$\tilde{V}=\{(a+bx)e+(a+b(x+1))\sigma~|~a,b\in k\}$$
Now we take the following $k$-basis of $\tilde{V}$:
$$\{u=e+\sigma,v=xe+(x+1)\sigma\}$$
If we write the $k$-endomorphisms of $\tilde{V}$ as matrices with respect
to the basis $\{u,v\}$, then we find the following grading on
on $B={\rm M}_2(k)$:
\begin{equation}\eqlabel{3.8.1}
B_e=k\pmatrix{1&0\cr 0&1\cr}\oplus k\pmatrix{0&\alpha\cr 1&1\cr}~~;~~
B_\sigma=k\pmatrix{1&1\cr 0&1\cr}\oplus k\pmatrix{0&\alpha\cr 1 &0\cr}
\end{equation}
As a consequence, we have now the list of all (different) isomorphism
types of $C_2$-graded algebra structures on ${\rm M}_2(k)$, with $k$ a
field of characteristic $2$, namely
\begin{itemize}
\item the trivial grading: all matrices have degree $e$;
\item the good grading ${\rm deg}(e_{11})={\rm deg}(e_{22})=e$,
${\rm deg}(e_{12})={\rm deg}(e_{21})=\sigma$;
\item the bad gradings from \eqref{3.8.1}, where $\alpha$ ranges over
a system of representatives of the $\{x^2-x~|~x\in k\}$ cosets of
$k$ different from $\{x^2-x~|~x\in k\}$.
\end{itemize}
Using direct computations, these gradings have been described before
in \cite{dinr}, and classified up to isomorphism in \cite{Boboc}.
\end{example}

\end{document}